\documentclass[12pt,reqno]{amsart}
\usepackage{graphicx}
\usepackage{times}
\usepackage{mathptmx}

\newtheorem{thm}{Theorem}[section]
\newtheorem{prop}[thm]{Proposition}
\newtheorem{lemma}[thm]{Lemma}

\newtheorem{cor}[thm]{Corollary}
\newtheorem{open}[thm]{Open Problem}

\theoremstyle{definition}
\newtheorem{dfn}[thm]{Definition}
\newtheorem{example}[thm]{Example}
\newtheorem{remark}[thm]{Remark}

\newcommand{\ER}{Erd\H{o}s-R\'enyi}
\newcommand{\bd}[1]{\mathbf{#1}}
\newcommand{\bdG}{\ensuremath{\bd G}}
\newcommand{\bdv}{\ensuremath{\bd v}}
\newcommand{\GG}{\mathcal{G}}
\newcommand{\GGn}{\ensuremath{\GG_n}}
\newcommand{\RR}{\mathbb{R}}
\newcommand{\XX}{\mathcal{X}}
\newcommand{\YY}{\mathcal{Y}}
\newcommand{\SSS}{\mathcal{S}}

\newcommand{\Gnp}{\ensuremath{\bd G(n,p)}}

\newcommand{\TV}{d_{\text{\rm\scriptsize TV}}}

\DeclareMathOperator{\EE}{E}

\newcommand{\AB}{B} 
\newcommand{\hatg}{\hat g}

\begin{document}

\title{On Vertex, Edge, and Vertex-Edge Random Graphs}

\author[Beer]{Elizabeth Beer}

\address{Center for Computing Sciences, 17100 Science Drive, Bowie, Maryland 20715-4300 USA}

\thanks{Elizabeth Beer's research on this paper, begun while she was a Ph.D.\ student at The Johns Hopkins University,
  was supported by a National Defense Science and Engineering Graduate Fellowship.  Research for James
  Fill is supported by The Johns Hopkins University's Acheson
  J.~Duncan Fund for the Advancement of Research in Statistics.}
  
\email{libby.beer@gmail.com}

 \author[Fill]{James Allen Fill}

\address{Department of Applied Mathematics and Statistics, The Johns
  Hopkins University, 3400 N.~Charles Street, Baltimore, Maryland 21218-2682 USA}

\email{jimfill@jhu.edu, ers@jhu.edu}

\author[Janson]{Svante Janson} 

\address{Department of Mathematics, Uppsala University, P.O. Box 480,
  SE-751 06 Uppsala, Sweden}

\email{svante.janson@math.uu.se} 

\author[Scheinerman]{Edward R. Scheinerman}



\begin{abstract}
  We consider three classes of random graphs: edge random graphs,
  vertex random graphs, and vertex-edge random graphs. Edge random
  graphs are \ER\ random graphs \cite{ER1,ER2}, vertex random graphs are
  generalizations of geometric random graphs \cite{penrose}, and
  vertex-edge random graphs generalize both. The names of
  these three types of random graphs describe where the randomness in
  the models lies: in the edges, in the vertices, or in both. We show
  that vertex-edge random graphs, ostensibly the most general of the
  three models, can be approximated arbitrarily closely by vertex
  random graphs, but that the two categories are distinct.
\end{abstract}

\date{October 12, 2010}

\maketitle


\section{Introduction}

The classic random graphs are those of Erd\H{o}s and R\'enyi
\cite{ER1,ER2}. In their model, each edge is chosen independently of
every other.  The randomness inhabits the edges; vertices simply serve
as placeholders to which random edges attach.

Since the introduction of \ER\ random graphs, many other models of
random graphs have been developed. For example, \emph{random geometric
  graphs} are formed by randomly assigning points in a Euclidean space
to vertices and then adding edges deterministically between vertices
when the distance between their assigned points is below a fixed
threshold; see \cite{penrose} for an overview. For these random
graphs, the randomness inhabits the vertices and the edges reflect
relations between the randomly chosen structures assigned to them.

Finally, there is a class of random graphs in which randomness is
imbued both upon the vertices and upon the edges. For example, in
latent position models of social networks, we imagine each vertex as
assigned to a random position in a metric ``social'' space. Then,
given the positions, vertices whose points are near each other are
more likely to be adjacent. See, for example, 
\cite{beer-thesis, HRH, rdpg, LSz, nickel-thesis}.  Such random graphs are, roughly
speaking, a hybrid of \ER\ and geometric graphs. 

We call these three categories, respectively, edge random, vertex
random, and vertex-edge random graphs. From their formal definitions
in Section~\ref{sect:rg}, it follows immediately that vertex random
and edge random graphs are instances of the more generous vertex-edge
random graph models. But is the vertex-edge random graph category
strictly more encompassing? We observe in Section~\ref{sect:approx}
that a vertex-edge random graph can be approximated arbitrarily
closely by a vertex random graph. Is it possible these two categories
are, in fact, the same? The answer is no, and this is presented in
Section~\ref{sect:not-the-same}. Our discussion closes in
Section~\ref{sect:open} with some open problems.

Nowadays, in most papers on random graphs, for each value of~$n$ a distribution is placed on the collection of
$n$-vertex graphs and asymptotics as $n \to \infty$ are studied.  We emphasize that in this paper, by contrast,
the focus is on what kinds of distributions arise in certain ways for a single arbitrary but fixed value of~$n$.

\section{Random Graphs}
\label{sect:rg}

For a positive integer~$n$, let $[n] = \{1, 2, \ldots, n\}$ and let
\GGn\ denote the set of all simple graphs $G = (V, E)$ with vertex set
$V = [n]$.  (A simple graph is an undirected graph with no loops and
no parallel edges.)  We often abbreviate the edge (unordered pair)
$\{i, j\}$ as $i j$ or write $i \sim j$ and say that~$i$ and~$j$ are
adjacent.

When we make use of probability spaces, we omit discussion of
measurability when it is safe to do so.  For example, when the sample
space is finite it goes without saying that the corresponding
$\sigma$-field is the total $\sigma$-field, that is, that all subsets
of the sample space are taken to be measurable.

\begin{dfn}[Random graph]
  \label{def:rg}
  A \emph{random graph} is a probability space of the form $\bd G =
  (\GGn, P)$ where~$n$ is a positive integer and~$P$ is a probability
  measure defined on \GGn.
\end{dfn}

In actuality, we \emph{should} define a random graph as a graph-valued
random variable, that is, as a measurable mapping from a probability
space into \GGn. However, the distribution of such a random object is
a probability measure on \GGn\ and is all that is of interest in this
paper, so the abuse of terminology in Definition~\ref{def:rg} serves
our purposes.

\begin{example}[\ER\ random graphs]
  \label{ex:ERG}
  A simple random graph is the \ER\ random graph in the case
  $p=\frac12$. This is the random graph $\bd G=(\GGn,P)$ where
  $$
  P(G) := 2^{-\binom n2}, \quad G \in \GGn.
  $$
  [Here and throughout we abbreviate $P(\{G\})$ as $P(G)$; this will
  cause no confusion.]  More generally, an \ER\ random graph is a
  random graph $\Gnp = (\GGn,P)$ where $p \in [0, 1]$ and
  $$
  P(G) := p^{|E(G)|} (1-p)^{\binom n2 - |E(G)|}, \quad G \in \GGn.
  $$
This means that the $\binom n2$ potential edges appear independently
of each other, each with probability~$p$.

This random graph model was first introduced by Gilbert~\cite{Gilbert}.
Erd\H{o}s and R\'enyi~\cite{ER1, ER2}, who started the systematic study of random
graphs, actually considered a closely related model with a fixed number of edges.
However, it is now common to call both models \ER\ random graphs.
\end{example}

\begin{example} [Single coin-flip random graphs] 
  \label{ex:coin-flip} 
  Another simple family of random graphs is one we call the
  \emph{single coin-flip} family.  Here $\bd G =(\GGn, P)$ where
  $p \in [0, 1]$ and
  $$
  P(G) :=
  \begin{cases}
    p & \text{if $G=K_n$,} \\
    1 - p & \text{if $G=\overline{K_n}$,} \\
    0 & \text{otherwise.}
  \end{cases}
  $$
As in the preceding example, each edge appears with probability~$p$;
but now all edges appear or none do.
\end{example}

In the successive subsections we specify our definitions of
\emph{edge}, \emph{vertex}, and \emph{vertex-edge} random graphs.

\subsection{Edge random graph}

In this paper, by an edge random graph (abbreviated ERG in the sequel)
we simply mean a classical \ER\ random graph.

\begin{dfn}[Edge random graph]
  An \emph{edge random graph} is an \ER\ random graph $\Gnp$.
\end{dfn}

We shall also make use of the following generalization that allows
variability in the edge-probabilities.

\begin{dfn}[Generalized edge random graph]
\label{def:GERG}
A generalized edge random graph (GERG) is a random graph for which the
events that individual vertex-pairs are joined by edges are mutually
independent but do not necessarily have the same probability.  Thus
to each pair $\{i, j\}$ of distinct vertices we associate a
probability $\bd p(i,j)$ and include
the edge $ij$ with probability $\bd p(i,j)$;
edge random graphs are
the special case where~$\bd p$ is constant.

Formally, a GERG can be described as follows.  Let~$n$ 
be a positive integer and let $\bd p: [n] \times [n] \to [0,1]$
be a symmetric function.  The \emph{generalized edge random graph} 
$\bd G( n,\bd p)$ is the probability space $(\GGn, P)$ with
  $$
  P(G) := \prod_{\substack{i<j \\ ij \in E(G)}} \bd p(i,j) \times
  \prod_{\substack{i<j \\ ij\notin E(G)}} [1-\bd p(i,j)].
  $$
\end{dfn}

We call the graphs in these two definitions (generalized) \emph{edge}
random graphs because all of the randomness inhabits the (potential)
edges.  The inclusion of ERGs in GERGs is strict, as easily constructed 
examples show.

GERGs have appeared previously
in the literature, e.g.\ in \cite{Alon};
see also the next example and Definition~\ref{def:rxphi} below.

\begin{example}[Stochastic blockmodel random graphs]
  \label{ex:stoch-block}
  A stochastic blockmodel random graph is 
  a GERG in which the vertex set is partitioned into blocks $B_1, B_2, \dots,
  B_b$ and the probability that vertices~$i$ and~$j$ are adjacent depends
  only on the blocks in which~$i$ and~$j$ reside.
  
  A simple example is a random bipartite graph defined by partitioning
  the vertex set into $B_1$ and $B_2$ and taking $\bd p(i,j)=0$ if
  $i,j\in B_1$ or $i,j\in B_2$, while $\bd p(i,j)=p$ (for some given
  $p$) if $i\in B_1$ and $j\in B_2$ or vice versa.

  The concept of blockmodel is interesting and useful when $b$
  remains fixed and $n \to \infty$.
  Asymptotics of blockmodel random graphs
  have been considered, for example, by S\"oderberg
  \cite{Sod1}. 
  (He also considers the version where the partitioning is random,
  constructed by independent random choices of a type in $\{1,...,b\}$
  for each vertex; see Example \ref{ex:VERG-finite}.)

  Recall, however, that in this paper we hold~$n$ fixed and note that in fact
  every GERG can be represented as a blockmodel by taking each
  block to be a singleton.

\end{example}

A salient feature of Example~\ref{ex:stoch-block} is that vertex labels matter. Intuitively, we
may expect that if all isomorphic graphs are treated ``the same'' by a
GERG, then it is an ERG.  We proceed to formalize this correct
intuition, omitting the simple proof of Proposition~\ref{prop:edge+homo}.

\begin{dfn}[Isomorphism invariance]
  Let $\bdG= (\GGn,P)$ be a random graph. We say that \bdG\ is
  \emph{isomorphism-invariant} if for all $G, H \in \GGn$ we have
  $P(G) = P(H)$ whenever~$G$ and~$H$ are isomorphic.
\end{dfn}
\begin{prop}
  \label{prop:edge+homo}
  Let \bdG\ be an isomorphism-invariant generalized edge random graph. 
  Then $\bdG = \Gnp$ for some $n, p$. That is, \bdG\ is an edge random graph.~\qed 
\end{prop}

\subsection{Vertex random graph}

The concept of a vertex random graph (abbreviated VRG) is motivated by
the idea of a random intersection graph.  One imagines a
universe~$\SSS$ of geometric objects.  A random $\SSS$-graph $G \in
\GGn$\ is created by choosing~$n$ members of~$\SSS$ independently at
random\footnote{Of course, some probability distribution must be
  associated with~$\SSS$.}, say $S_1, \ldots, S_n$, and then declaring
distinct vertices~$i$ and~$j$ to be adjacent if and only if $S_i \cap
S_j \not= \emptyset$.  For example, when~$\SSS$ is the set of real
intervals, one obtains a random interval graph 
\cite{DHJinterval, random-intervals-monthly,rand-interval,rigs-evolution}; 
see Example~\ref{E:randint} for more.
In
\cite{fss,kss,singer-thesis} one takes~$\SSS$ to consist of discrete
(finite) sets. Random chordal graphs can be defined by selecting
random subtrees of a tree~\cite{rand-chordal}.

Notice that for these random graphs, all the randomness lies in the
structures attached to the vertices; once these random structures have
been assigned to the vertices, the edges are \emph{determined}.  In
Definition~\ref{def:VRG} we generalize the idea of a random
intersection graph to other vertex-based representations of graphs;
see~\cite{spinrad}.


\begin{dfn}[$({\bd x}, \phi)$-graph]
  \label{def:xphi}
  Let~$n$ be a positive integer, $\XX$ a set, ${\bd x} = (x_1, \dots,
  x_n)$ a function from $[n]$ into~$\XX$, and $\phi:\XX \times \XX \to
  \{0,1\}$ a symmetric function.  Then the \emph{$({\bd x},
    \phi)$-graph}, denoted $\bdG(\bd x, \phi)$, is defined to be the
  graph with vertex set $[n]$ such that for all $i, j \in [n]$ with $i
  \neq j$ we have
  $$
  i j \in E\quad\mbox{if and only if}\quad\phi(x_i, x_j) = 1.
  $$
\end{dfn}

Of course, every graph $G = (V, E)$ with $V = [n]$ is an $({\bd x},
\phi)$-graph for some choice of $\XX$, ${\bd x}$, and~$\phi$; one need
only take~${\bd x}$ to be the identity function on $\XX := [n]$ and
define
$$
\phi(i, j) := {\bd 1}(i j \in E) = 
\begin{cases}
    1 & \text{if $i j \in E$} \\
    0 & \text{otherwise.}
\end{cases}
$$
It is also clear that this representation of~$G$ as an $({\bd x},
\phi)$-graph is far from unique.  The notion of $({\bd x},
\phi)$-graph becomes more interesting when one or more of $\XX$, ${\bd
  x}$, and~$\phi$ are specified.

\begin{example}[Interval graphs]
  \label{ex:random-intervals}
  Take~$\XX$ to be the set of all real intervals and define
  \begin{equation}
  \label{randintphi}
  \phi(J, J') := 
  \begin{cases}
    1 & \text{if $J \cap J' \not=\emptyset$} \\
    0 & \text{otherwise.}
  \end{cases}
  \end{equation}
  In this case, an $({\bd x}, \phi)$-graph is exactly an interval graph.
\end{example}

\begin{dfn}[Vertex random graph]
  \label{def:VRG}
To construct a vertex random graph (abbreviated VRG),
we imbue~$\XX$ with a probability measure~$\mu$ and sample~$n$
elements of~$\XX$ independently at random to get~${\bd x}$,
and then we build the $({\bd x}, \phi)$-graph.

Formally, let~$n$ be a positive integer, $(\XX,\mu)$ a probability space, and
  $\phi:\XX \times \XX \to \{0,1\}$ a symmetric function.  The
  \emph{vertex random graph} $\bd G(n, \XX, \mu, \phi)$ is the random
  graph $(\GGn, P)$ with
  $$
  P(G) := \int\!\bd1\{\bdG(\bd x, \phi) = G\}\,{\mathbf \mu}(d {\bd
    x}),\quad G \in \GGn,
  $$
  where ${\mathbf \mu}(d {\bd x})$ is shorthand for the product
  integrator $\mu^n(d {\bd x}) = \mu(dx_1) \dots \mu(dx_n)$ on
  $\XX^n$.
\end{dfn}

Note that $\bdG(\cdot, \phi)$ is a graph-valued random variable
defined on $\XX^n$.  The probability assigned by the vertex random
graph to $G \in \GGn$ is simply the probability that this random
variable takes the value~$G$.

\begin{example}[Random interval graphs]
\label{E:randint}
Let~$\XX$ be the set of real intervals as in Example~\ref{ex:random-intervals},
let~$\phi$ be as in~\eqref{randintphi}, and let~$\mu$ be a probability measure
on~$\XX$.  This yields a VRG that is a random interval graph.
\end{example}

\begin{example}[Random threshold graphs]
\label{E:thresh}
Let~$\XX = [0, 1]$, let $\mu$ be Lebesgue measure, and let~$\phi$ be the
indicator of a given up-set in the usual (coordinatewise) partial order on $\XX \times \XX$.
This yields a VRG that is a random threshold graph; see~\cite{DHJthreshold}.
\end{example}

\begin{example}[Random geometric graphs]
  Random geometric graphs are studied extensively in~\cite{penrose}.
  Such random graphs are created by choosing $n$ i.i.d.\ (independent
  and identically distributed) points from some probability
  distribution on $\RR^k$.  Then, two vertices are joined by an edge
  exactly when they lie within a certain distance, $t$, of each other.

  Expressed in our notation, we let $(\XX,d)$ be a metric space
  equipped with a probability measure~$\mu$ and let $t>0$ (a
  threshold). For points $x,y\in\XX$ define
  $$
  \phi(x,y) := \bd1 \left\{ d(x,y) \le t \right\}.
  $$
  That is, two vertices are adjacent exactly when the distance between
  their corresponding randomly chosen points is sufficiently small.
\end{example}

Because the~$n$ vertices in a vertex random graph are drawn i.i.d.\
from $(\XX,\mu)$, it is easy to see that the random graph is
isomorphism-invariant.

\begin{prop}
  \label{prop:vertex-homo}
  Every vertex random graph is isomorphism-invariant.\qed
\end{prop}

\subsection{Vertex-edge random graphs}

A generalization both of vertex random graphs and of edge random
graphs are the \emph{vertex-edge} random graphs (abbreviated VERGs) of
Definition~\ref{def:VERG}.  First we generalize
Definition~\ref{def:xphi} to allow edge probabilities other than~$0$
and~$1$.

\begin{dfn}[Random $({\bd x}, \phi)$-graph]
  \label{def:rxphi}
Given a positive integer $n \geq 1$, a set~$\XX$,
and a function $\phi:\XX \times \XX \to [0, 1]$,
we assign to each $i \in [n]$ a deterministically chosen
object $x_i \in \XX$.  Then, for each pair $\{i, j\}$ of vertices,
independently of all other pairs, the edge $i j$ is included
in the random $({\bd x}, \phi)$-graph with probability $\phi(x_i, x_j)$.  

Formally, let ${\bd x} = (x_1, \dots, x_n)$ be a given function 
from $[n]$ into~$\XX$.  Then the 
\emph{random $({\bd x}, \phi)$-graph}, denoted $\bdG(\bd x, \phi)$, 
is defined to be the
random graph $(\GGn, P_{\bd x})$ for which the probability of $G \in
\GGn$ is given by
  $$
  P_{\bd x}(G) := \prod_{i < j,\ i \sim j} \phi(x_i, x_j) \times
  \prod_{i < j,\ i \not\sim j} [1-\phi(x_i, x_j)].
  $$
\end{dfn}

Notice that $\bdG(\bd x, \phi)$ is simply the generalized
edge random graph $\bd G(n, \bd p)$ where $\bd p(i, j) := \phi(x_i,
x_j)$ (recall Definition~\ref{def:GERG}).

\begin{dfn}[Vertex-edge random graph]
  \label{def:VERG}
Let~$n$ be a positive integer, $(\XX,\mu)$ a probability space, and
$\phi:\XX \times \XX \to [0,1]$ a symmetric function.  
In words, a vertex-edge random graph is generated like this:
First a list of random elements is drawn i.i.d.,\ with distribution~$\mu$,
from~$\XX$; call the
list $\bd X = (X_1, \ldots, X_n)$.  Then, conditionally given $\bd X$,
independently for each pair of distinct vertices~$i$ and~$j$ we
include the edge $i j$ with probability $\phi(X_i, X_j)$.

Formally, the \emph{vertex-edge random graph} $\bd G(n, \XX, \mu, \phi)$ is the
random graph $(\GGn, P)$ with
  $$
  P(G) := \int\!P_{\bd x}(G)\,{\mathbf \mu}(d {\bd x})
  $$
where the integration notation is as in Definition~\ref{def:VRG} and
$P_{\bd x}$ is the probability measure for the random 
$({\bd x}, \phi)$-graph $\bd G(\bd x, \phi)$ of Definition~\ref{def:rxphi}.
\end{dfn}

Note that a VRG is the special case of a VERG with~$\phi$ taking values
in $\{0, 1\}$. 

It can be shown \cite{SJ210} that every VERG can be constructed with
the standard choice $\XX=[0,1]$ and $\mu =$ Lebesgue measure. However,
other choices are often convenient in specific situations.

We note in passing that one could generalize the notions of VRG and
VERG in the same way that edge random graphs (ERGs) were generalized
in Definition~\ref{def:GERG}, by allowing different functions $\phi_{i j}$ for
different vertex pairs $\{i, j\}$.  But while the notion of generalized ERG
was relevant to the definition of a VERG (recall the sentence
preceding Definition~\ref{def:VERG}), we neither study nor employ
generalized VRGs and VERGs in this paper.

Asymptotic properties (as $n\to\infty$) of random $({\bd x},
\phi)$-graphs and VERGs have been studied by several authors:
see, e.g.,\ \cite{SJ178} and the references therein.  VERGs are also important
in the theory of \emph{graph limits}; see for example
\cite{BCL1,SJ209,LSz}.

\begin{example}[Finite-type VERG]
  \label{ex:VERG-finite}
  In the special case when $\XX$ is finite, $\XX=\{1,\ldots,b\}$ say,
  we thus randomly and independently choose a type in $\{1,\ldots,b\}$
  for each vertex, with a given distribution $\mu$; we can regard this
  as a random partition of the vertex set into blocks $B_1,\ldots,B_b$
  (possibly empty, and with sizes governed by a multinomial
  distribution). A VERG with $\XX$ finite can thus be regarded as a
  stochastic blockmodel graph with multinomial random blocks; cf.\
  Example~\ref{ex:stoch-block}.  Such finite-type VERGs have been
  considered by S\"oderberg \cite{Sod1,Sod3,Sod4,Sod2}.
\end{example}

\begin{example} [Random dot product graphs]
  \label{ex:rpdg}
  In \cite{rdpg,nickel-thesis} random graphs are generated by the
  following two-step process.  First, $n$ vectors (representing $n$
  vertices) $\bdv_1,\ldots,\bdv_n$ are chosen i.i.d.\ according to
  some probability distribution on $\RR^k$.  With this choice in
  place, distinct vertices~$i$ and~$j$ are made adjacent with
  probability $\bdv_i \cdot \bdv_j$. All pairs are considered
  (conditionally) independently. Care is taken so that the
  distribution on $\RR^k$ satisfies
  $$
  P\bigl( \bd v_i \cdot \bd v_j \notin [0,1] \bigr) = 0.
  $$

  Random dot product graphs are vertex-edge random graphs with
  $\XX=\RR^k$ and $\phi(\bd v,\bd w) = \bd v \cdot \bd w$. 
\end{example}

As with vertex random graphs, all vertices are treated ``the same'' in
the construction of a vertex-edge random graph.

\begin{prop}
  Every vertex-edge random graph is isomorphism-invariant.\qed
\end{prop}

Note that we use the notation $\bd G(n,\XX,\mu,\phi)$ for both VRGs
and VERGs.  This is entirely justified because~$\phi$ takes values in
in $\{0, 1\}$ for VRGs and in $[0,1]$ for VERGs.  If perchance
the~$\phi$ function for a VERG takes only the values~$0$ and~$1$, then
the two notions coincide.  Hence we have part~(b) of the following
proposition; part~(a) is equally obvious.

\begin{prop}
  \label{prop:VERGinclusions}
  \
  \begin{itemize}
  \item [(a)]
  Every edge random graph is a vertex-edge random graph.

  \item [(b)] Every vertex random graph is a vertex-edge random graph.
  \qed
  \end{itemize}
\end{prop}

However, not all generalized edge random graphs are vertex-edge random
graphs, as simple counterexamples show.

We now ask whether the converses to the statements in
Proposition~\ref{prop:VERGinclusions} are true.
The converse to Proposition~\ref{prop:VERGinclusions}(a) is false.
Indeed, It is easy to find examples of VERGs that are not ERGs:
\begin{example}
We present one small class of examples of VERGs that are even VRGs,
but not ERGs.  Consider
  \emph{random interval graphs} \cite{DHJinterval, random-intervals-monthly,
    rand-interval} $\bd G(n,\XX,\mu,\phi)$ with $n \geq 3$, $\XX$
  and~$\phi$ as in Example~\ref{ex:random-intervals}, and (for $i \in
  [n]$) the random interval $J_i$ corresponding to vertex~$i$
  constructed as $[X_i, Y_i]$ or $[Y_i, X_i]$, whichever is nonempty,
  where $X_1, Y_1, \dots, X_n, Y_n$ are i.i.d.\ uniform$[0, 1]$ random
  variables.  From an elementary calculation, independent of~$n$, one
  finds that the events $\{1 \sim 2\}$ and $\{1 \sim 3\}$ are not
  independent.
\end{example}

The main result of this paper (Theorem~\ref{thm:not-the-same}; see
also the stronger Theorem~\ref{thm:ERG-not-VRG}) is that the converse
to Proposition~\ref{prop:VERGinclusions}(b) is also false.  The class
of vertex random graphs does not contain the class of vertex-edge
random graphs; however, as shown in the next section, every
vertex-edge random graph can be approximated arbitrarily closely by a
vertex random graph.

\begin{figure}[ht]
  \centering
  \includegraphics[width=\textwidth]{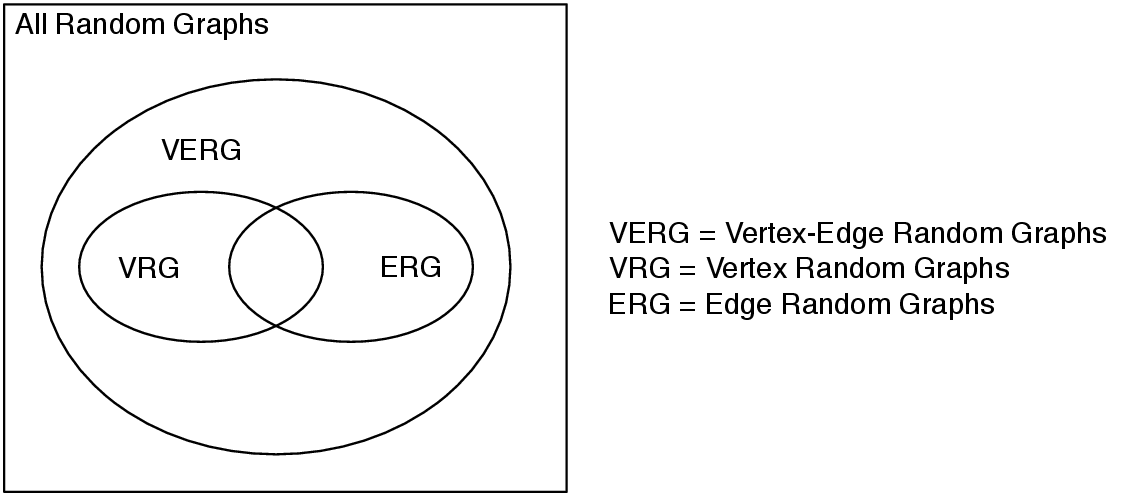}
  \caption{Venn diagram of random graph classes. The results of this
    paper show that all five regions in the diagram are nonempty.}
  \label{fig:venn}
\end{figure}

An overview of the inclusions of these various categories is presented
in Figure~\ref{fig:venn}.

\section{Approximation}
\label{sect:approx}

The goal of this section is to show that every vertex-edge random
graph can be closely approximated by a vertex random graph. Our notion
of approximation is based on total variation distance.  (This choice
is not important. We consider a fixed $n$, and the space of
probability measures on $\GGn$ is a finite-dimensional simplex, and
thus compact. Hence any continuous metric on the probability measures
on $\GGn$ is equivalent to the total variation distance, and can be
used in Theorem \ref{thm:approx}.)

\begin{dfn}[Total variation distance]
  \label{dfn:TV}
  Let $\bd G_1=(\GGn,P_1)$ and $\bd G_2=(\GGn,P_2)$ be random graphs
  on $n$ vertices. We define the \emph{total variation distance}
  between $\bd G_1$ and $\bd G_2$ to be 
  $$
  \TV(\bd G_1, \bd G_2) = \frac12\sum_{G \in \GGn} \left| P_1(G) -
    P_2(G) \right|.
  $$
\end{dfn}

Total variation distance can be reexpressed in terms of the maximum
discrepancy of the probability of events.


\begin{prop} \label{prop:TV} 
  Let $\bd G_1=(\GGn,P_1)$ and $\bd G_2=(\GGn,P_2)$ be random graphs
  on $n$ vertices. Then
  $$
  \TV(\bd G_1, \bd G_2) = \max_{\AB \subseteq \GGn} \left| P_1(\AB) -
    P_2(\AB) \right|. \qed
  $$
\end{prop}

\begin{thm}
  \label{thm:approx}
  Let $\bdG$ be a vertex-edge random graph and let
  $\varepsilon>0$. There exists a vertex random graph $\widehat\bdG$ with
  $\TV(\bdG, \widehat\bdG) < \varepsilon$. 
\end{thm}

We  use the following simple birthday-problem subadditivity upper bound.
Let $M$ be a positive integer.

\begin{lemma}
  \label{lem:no-repeat}
  Let $\bd A = (A_1, A_2, \ldots, A_n)$ be a random sequence of integers
  with each $A_i$ chosen independently and uniformly from $[M]$. Then
  $$
  P\left\{ \text{\rm $\bd A$ has a repetition} \right\}
  \leq \frac{n^2}{2M}.~\qed
  $$
\end{lemma}
\smallskip

\begin{proof}[Proof of Theorem~\ref{thm:approx}]
  Let~\bdG\ be a vertex-edge random graph on $n$ vertices and let
  $\varepsilon > 0$. Let~$M$ be a large positive integer. (We postpone
  our discussion of just how large to take~$M$ until needed.) 

  The vertex-edge random graph \bdG\ can be written $\bdG =
  \bdG(n,\XX,\mu,\phi)$ for some set $\XX$ and mapping
  $\phi:\XX\times\XX\to[0,1]$. 

  We construct a vertex random graph $\widehat\bdG =
  \bdG(n,\YY,\nu,\psi)$ as follows.  Let $\YY := \XX \times [0,1]^{M}
  \times [M]$; that is, $\YY$ is the set of ordered triples $(x, f,
  a)$ where $x \in \XX$, $f\in[0,1]^{M}$, and $a \in [M]$. We
  endow~$\YY$ with the product measure of its factors; that is, we
  independently pick $x \in \XX$ according to $\mu$, a function
  $f\in[0,1]^{[M]}$ uniformly, and $a \in [M]$ uniformly. We denote
  this measure by~$\nu$.

  We denote the components of the vector $f\in[0,1]^M$ by
  $f(1),\dots,f(M)$, thus regarding $f$ as a random function
  from~$[M]$ into $[0,1]$.  Note that for a random $f\in[0,1]^M$, the
  components $f(1),\dots,f(M)$ are i.i.d.\ random numbers with a
 uniform$[0, 1]$ distribution.

 Next we define $\psi$. Let $y_1,y_2 \in \YY$ where
 $y_i=(x_i, f_i, a_i)$ (for $i=1,2$). Let
  $$
  \psi(y_1,y_2) =
  \begin{cases}
    1 & \text{if $a_1 < a_2$ and $\phi(x_1,x_2)\ge f_1(a_2)$,} \\
    1 & \text{if $a_2 < a_1$ and $\phi(x_1,x_2)\ge f_2(a_1)$,} \\
    0 & \text{otherwise.}
  \end{cases}
  $$
  Note that~$\psi$ maps $\YY \times \YY$ into $\{0,1\}$ and is symmetric in its
  arguments. Therefore $\widehat\bdG$ is a vertex
  random graph. 
  \medbreak

  We now show that $\TV(\bdG, \widehat\bdG)$ can be made arbitrarily small
  by taking~$M$ sufficiently large.

  Let $\AB\subseteq \GG_n$. Recall that
  \begin{align*}
    P(\AB) &=
    \int\!P_{\bd x}(\AB)\,{\mathbf \mu}(d {\bd x}), \\
    \widehat{P}(\AB) &= \int\!\bd1\{\bdG(\bd y, \psi) \in \AB\}
    \,{\mathbf \nu}(d {\bd y}) = \Pr\{\bd G(\bd Y, \psi) \in \AB\},
  \end{align*}
  where in the last expression the~$n$ random variables comprising
  ${\bd Y} = (Y_1, \dots, Y_n)$ are independently 
  chosen from~$\YY$, each according to the distribution~$\nu$.

  As each $Y_i$ is of the form $(X_i, F_i, A_i)$ we break up the
  integral for $\widehat{P}(\AB)$ based on whether or not the
  $a$-values of the $Y$s are repetition free and apply
  Lemma~\ref{lem:no-repeat}:
  \begin{equation}
    \label{q1}
    \begin{split}
      \widehat P(\AB) &= \Pr\{\bd G(\bd Y, \psi) \in \AB \mid \bd A \text{ is
        repetition free}\} \Pr\{\bd A \text{ is repetition free}\}
      \\
      &\qquad+ \Pr\{\bd G(\bd Y, \psi) \in \AB \mid \bd A \text{ is not
        repetition free}\} \Pr\{\bd A \text{ is not repetition free}\} \\
      &=
      \Pr\{\bd G(\bd Y, \psi) \in \AB \mid \bd A \text{ is
        repetition free}\} + \delta
    \end{split}
  \end{equation}
  where $|\delta| \leq n^2 / (2 M)$.

  Now, for any repetition-free $\bd a$, the events $\{i \sim j\text{
    in $\bd G(\bd Y, \psi)$}\}$ are conditionally independent
  given~$\bd X$ and given $\bd A = \bd a$, with
  \begin{align*}
    \Pr\{i \sim j\mbox{\ in $\bd G(\bd Y, \psi)$} \mid \bd X,\ \bd A = \bd a\} &= 
    \begin{cases}
      \Pr\{\phi(X_i, X_j) \geq F_i(a_j) \mid X_i, X_j\} & \text{if\ }a_i < a_j \\
      \Pr\{\phi(X_i, X_j) \geq F_j(a_i) \mid X_i, X_j\} & \text{if\ }a_j < a_i 
    \end{cases} \\
    &= \phi(X_i, X_j).
  \end{align*}
  Thus, for any repetition-free $\bd a$, 
  $$
  \Pr\{\bd G(\bd Y, \psi) \in \AB \mid \bd X,\ \bd A = \bd a\}  
  $$
  equals
  $$
  \sum_{G\in\AB}\; \left( \prod_{i < j,\ ij\in E(G)} \phi(X_i,X_j)
    \times \prod_{i < j,\ ij \notin E(G)} \left[ 1 - \phi(X_i,X_j)
    \right] \right) = P_{\bd X} (\AB).
  $$
  Removing the conditioning on~$\bd X$ and~$\bd A$, \eqref{q1} thus
  implies
  $$
  \widehat P(\AB) = P(\AB) +  \delta,
  $$
  and so $|P(\AB) - \widehat P(\AB)| \leq n^2 / M$ for all $\AB
  \subseteq \GGn$. Equivalently, $\TV(\bdG, \widehat\bdG) \le n^2/M$.
  Thus we need only choose $M > n^2 / \varepsilon$.
\end{proof}

\section{Not all vertex-edge random graphs are vertex random graphs}
\label{sect:not-the-same}

In Section~\ref{sect:approx} (Theorem~\ref{thm:approx}) it was shown
that every vertex-edge random graph can be approximated arbitrarily
closely by a vertex random graph.  This naturally raises the question
of whether every vertex-edge random graph \emph{is} a vertex random
graph.  We originally believed that some suitable ``$M = \infty$
modification'' of the proof of Theorem~\ref{thm:approx} would provide
a positive answer, but in fact the answer is no:

\begin{thm}
  \label{thm:not-the-same}
  Not all vertex-edge random graphs are vertex random graphs.
\end{thm}

This theorem is an immediate corollary of the following much stronger
result.  We say that an ERG $\Gnp$ is \emph{nontrivial} when $p \notin
\{0, 1\}$.

\begin{thm}
  \label{thm:ERG-not-VRG}
  If $n \geq 4$, no nontrivial \ER\ random graph is a vertex random
  graph.  In fact, an ERG $\Gnp$ with $n \geq 4$ is represented as a
  vertex-edge random graph $\bd G(n, \XX, \mu, \phi)$ if and only if
  $\phi(x, y) = p$ for $\mu$-almost every~$x$ and~$y$.
\end{thm}

The ``if'' part of Theorem~\ref{thm:ERG-not-VRG} is trivial (for any
value of~$n$), since $\phi(x, y) = p$ clearly gives a representation
(which we shall call the \emph{canonical} representation) of an ERG as
a VERG.

We establish a lemma before proceeding to the proof of the
nontrivial ``only if'' part of Theorem~\ref{thm:ERG-not-VRG}.  To set
up for the lemma, which relates an expected subgraph count to
the spectral decomposition of a certain integral operator, consider
any particular representation $\bd G(n, \XX, \mu, \phi)$ of a VERG.
Let~$T$ be the integral operator with kernel~$\phi$ on the space
$L(\XX, \mu)$ of $\mu$-integrable functions on~$\XX$:
\begin{equation}
  \label{T}
  (T g)(x) := \int\!\phi(x, y) g(y)\,\mu(dy) = \EE [\phi(x, X) g(X)]
\end{equation}
where~$\EE$ denotes expectation and~$X$ has distribution~$\mu$.
Since~$\phi$ is bounded and symmetric and~$\mu$ is a finite measure,
$T$ is a self-adjoint Hilbert--Schmidt operator.  Let the finite or
infinite sequence $\lambda_1, \lambda_2, \dots$ denote its eigenvalues
(with repetitions if any); note that these are all real.  Note also
that in the special case $\phi(x, y) \equiv p$ giving the canonical
representation of an ERG, we have $\lambda_1 = p$ and $\lambda_i = 0$
for $i \geq 2$.

Let $N_k$, $3 \leq k \leq n$, be the number of rooted $k$-cycles in
$\bd G(n, \XX, \mu, \phi)$, where a (not necessarily induced) rooted
cycle is a cycle with a designated start vertex (the root) and a start
direction.  In the following we write $n^{\underline{k}} := n (n - 1)
\cdots (n - k + 1)$ for the $k$th falling factorial power of~$n$.

\begin{lemma}
  \label{lem:rtdcyc}
  In a VERG, with the preceding notation, for $3\le k\le n$ we have
  $$
  \EE N_k = n^{\underline{k}} \sum_i \lambda^k_i.
  $$
\end{lemma}

\begin{proof}
  A rooted $k$-cycle is given by a sequence of~$k$ distinct vertices
  $v_1, \dots, v_k$ with edges $v_i v_{i + 1}$ ($i = 1, \dots, k - 1$)
  and $v_k v_1$.  Thus, with $\mbox{Tr}$ denoting trace,
  \begin{align*}
    \EE N_k &= 
    n^{\underline{k}} 
    \EE[\phi(X_1, X_2) \phi(X_2, X_3) \cdots \phi(X_k, X_1)] \\
    &= n^{\underline{k}} \int \cdots \int_{\XX^k} \phi(x_1, x_2)
    \phi(x_2, x_3) \cdots
    \phi(x_k, x_1)\,d\mu(x_1) \cdots d\mu(x_k) \\
    &= n^{\underline{k}} \mbox{Tr}\,T^k = n^{\underline{k}} \sum_i
    \lambda^k_i.  \qedhere
  \end{align*}
\end{proof}

In the special case $\phi(x, y) \equiv p$ of the canonical
representation of an ERG, Lemma~\ref{lem:rtdcyc} reduces to
\begin{equation}
  \label{ERGrc}
  \EE N_k = n^{\underline{k}} p^k, \qquad 3 \leq k \leq n,
\end{equation}
which is otherwise clear for an ERG.

Equipped with Lemma~\ref{lem:rtdcyc}, it is now
easy to prove Theorem~\ref{thm:ERG-not-VRG}.  

\medskip

\begin{proof}[Proof of Theorem~\ref{thm:ERG-not-VRG}]
In any VERG $\bd G(n, \XX, \mu, \phi)$, the average edge-prob\-a\-bil\-i\-ty~$\rho$ is given by
$$
\rho := \EE \phi(X_1, X_2) = \int\!\int\!\phi(x, y)\,\mu(dy)\,\mu(dx) = \langle T {\bf 1}, {\bf 1} \rangle \leq \lambda_1,
$$
where $\bf 1$ is the function with constant value~$1$ and~$\lambda_1$ is the largest eigenvalue of~$T$; hence
\begin{equation}
\label{ineqs}
\rho^4 \leq \lambda^4_1 \leq \sum_i \lambda^4_i = \frac{\EE N_4}{n^{\underline{4}}},
\end{equation}
where the equality here comes from Lemma \ref{lem:rtdcyc}.  If the VERG is an ERG $\Gnp$,
then $\rho = p$ and by combining~\eqref{ERGrc} and~\eqref{ineqs} we see that $p = \rho = \lambda_1$
and $\lambda_i = 0$ for $i \geq 2$; 
hence, $\phi(x, y) = p \psi(x) \psi(y)$ for $\mu$-almost every~$x$ and~$y$,
where~$\psi$ is a normalized eigenfunction of~$T$ corresponding to eigenvalue $\lambda_1 = p$.
But then
$$
p \int\!\psi^2(x)\,\mu(dx) = p = \int\!\int\!\phi(x, y)\,\mu(dy)\,\mu(dx) = p \left[ \int\!\psi(x)\,\mu(dx) \right]^2,
$$
and since there is equality in the Cauchy--Schwarz inequality for~$\psi$ we see that
$\phi(x, y) = p$ for $\mu$-almost every~$x$ and~$y$.  This establishes the ``only if''
assertion in Theorem~\ref{thm:ERG-not-VRG}; as already noted, the ``if'' assertion is trivial.
\end{proof}

Consider an ERG $\Gnp$.  If $n \geq 4$, Theorem~\ref{thm:ERG-not-VRG}
shows that $\Gnp$ is never a VRG if $p \notin \{0, 1\}$.  Curiously,
however, every $\Gnp$ with $n \leq 3$ is a VRG; in fact, the following
stronger result is true.

\begin{thm}
  \label{thm:VERG-is-VRG}
  Every vertex-edge random graph with $n \leq 3$ is a vertex random graph.
\end{thm}

\begin{proof}
  We seek to represent the given VERG $\bd G(n, \XX, \mu, \phi)$ as a
  VRG $\bd G(n, \YY, \nu, \psi)$, with $\psi$ taking values in $\{0,
  1\}$.  For $n = 1$ there is nothing to prove.  For $n = 2$, the only
  random graphs of any kind are ERGs $\Gnp$; one easily checks that
  $\YY$ = $\{0, 1\}$, $\nu({1}) = \sqrt{p} = 1 - \nu({0})$, and
  $\psi(y_1, y_2) = {\bf 1}(y_1 = y_2 = 1)$ represents $\Gnp$ as a
  VRG.

  Suppose now that $n = 3$.  The given VERG can be described as
  choosing $X_1, X_2, X_3$ i.i.d.\ from~$\mu$ and, independently,
  three independent uniform$[0, 1)$ random variables $U_{1 2}, U_{1
    3}, U_{2 3}$, and then including each edge $i j$ if and only if
  the corresponding $U_{i j}$ satisfies $U_{i j} \leq \phi(X_i, X_j)$.
  According to Lemma~\ref{lem:unif} to follow, we can obtain such
  $U_{i j}$'s by choosing independent uniform$[0, 1)$ random
  variables $U_1, U_2, U_3$ and setting $U_{i j} := U_i \oplus U_j$,
  where $\oplus$ denotes addition modulo~$1$.  It follows that the
  given VERG is also the VRG $\bd G(3, \YY, \nu, \psi)$, where $\YY :=
  \XX \times [0, 1)$, $\nu$ is the product of~$\mu$ and the
  uniform$[0, 1)$ distribution, and, with $y_i = (x_i, u_i)$,
  \begin{equation}
  \label{psi}
  \psi(y_1, y_2) = {\bf 1}(u_1 \oplus u_2 \leq \phi(x_1, x_2)).
  \end{equation}
\end{proof}
\smallskip

\begin{lemma}
  \label{lem:unif}
  If $U_1, U_2, U_3$ are independent uniform$[0, 1)$ random
  variables, then so are $U_1 \oplus U_2$, $U_1 \oplus U_3$, $U_2
  \oplus U_3$, where $\oplus$ denotes addition modulo~$1$.
\end{lemma}

\begin{proof}
  The following proof seems to be appreciably simpler than a
  change-of-variables proof.  For other proofs, see Remark
  \ref{Runiform} below.  Let $J := \{0, \dots, k - 1\}$.  First check
  that, for $k$ odd, the mapping
  $$
  (z_1, z_2, z_3) \mapsto (z_1 + z_2, z_1 + z_3, z_2 + z_3),
  $$
  from $J \times J \times J$ into $J \times J \times J$, with addition
  here modulo~$k$, is bijective.  Equivalently, if $U_1, U_2, U_3$ are
  iid uniform$[0, 1)$, then the joint distribution of
  \begin{eqnarray*}
  Z_{1 2}(k) &:=& \lfloor k U_1 \rfloor + \lfloor k U_2 \rfloor, \\
  Z_{1 3}(k) &:=& \lfloor k U_1 \rfloor + \lfloor k U_3 \rfloor, \\
  Z_{2 3}(k) &:=& \lfloor k U_2 \rfloor + \lfloor k U_3 \rfloor
  \end{eqnarray*}
  is the same as that of
  $$
  \lfloor k U_1 \rfloor,   \lfloor k U_2 \rfloor,   \lfloor k U_3 \rfloor .
  $$
  Dividing by~$k$ and letting $k \to \infty$ through odd values of~$k$
  gives the desired result.
\end{proof}

\begin{remark}
  \label{Rhyper}
  Theorem~\ref{thm:VERG-is-VRG} has an extension to hypergraphs.
  Define a VERHG (vertex-edge random hypergraph) on the vertices $\{1,
  \dots, n\}$ in similar fashion to VERGs, except that now each of
  the~$n$ possible hyperedges joins a subset of vertices of size $n -
  1$.  Define a VRHG (vertex random hypergraph) similarly.  Then
  VERHGs and VRHGs are the same, for each fixed $n$.  The key to the
  proof is the observation (extending the case $n = 3$ of
  Lemma~\ref{lem:unif}) that if $U_1, U_2, \dots U_n$ are i.i.d.\
  uniform$[0, 1)$, then the same is true (modulo~$1$) of $S - U_1, S
  - U_2, \dots, S - U_n$, where $S := U_1 + U_2 + \cdots + U_n$.  The
  observation can be established as in the proof of
  Lemma~\ref{lem:unif}, now by doing integer arithmetic modulo~$k$,
  where $n - 1$ and~$k$ are relatively prime, and passing to the limit
  as $k \to \infty$ through such values.  [For example, consider $k =
  m (n - 1) + 1$ and let $m \to \infty$.]
\end{remark}

\begin{remark}\label{Runiform}
  Consider again Lemma~\ref{lem:unif} and its extension in Remark
  \ref{Rhyper}.  Let~${\bf T} = {\bf R}/{\bf Z}$ denote the circle.
  We have shown that the mapping ${\bf u} \mapsto A {\bf u}$ preserves
  the uniform distribution on ${\bf T}^n$, where for example in the
  case $n = 3$ the matrix~$A$ is given by
  $$
  A = \left( 
    \begin{array}{ccc}
      1 & 1 & 0 \\
      1 & 0 & 1 \\
      0 & 1 & 1
    \end{array}
  \right).
  $$
  More generally, the mapping ${\bf u} \mapsto A {\bf u}$ preserves
  the uniform distribution on ${\bf T}^n$ whenever~$A$ is a
  nonsingular $n \times n$ matrix of integers.  Indeed, then $A: {\bf
    R}^n \to {\bf R}^n$ is surjective, so $A: {\bf T}^n \to {\bf T}^n$
  is surjective; and any homomorphism of a compact group (here ${\bf
    T}^n$) onto a compact group (here also ${\bf T}^n$) preserves the
  uniform distribution, i.e., the (normalized) Haar measure.  (This
  follows, e.g.,\ because the image measure is translation invariant.)
  This preservation can also be seen by Fourier analysis: For the
  i.i.d.\ uniform vector ${\bf U} = (U_1, \dots, U_n)$ and any integer
  vector ${\bf k} = (k_1, \dots, k_n) \neq {\bf 0}$,
  $$
  \EE \exp(2 \pi i {\bf k} \cdot A {\bf U}) = \EE \exp(2 \pi i A^T {\bf k} \cdot {\bf U}) = 0
  $$
  because $A^T {\bf k} \neq 0$.
\end{remark}

\begin{remark}
  \label{rmk:n=3}
  In this remark we (a)~give a spectral characterization of all
  representations of a three-vertex ERG ${\bd G(3, p)}$ as a VERG $\bd
  G(3, \XX, \mu, \phi)$ and (b)~briefly discuss the spectral
  decomposition of the ``addition modulo 1'' kernel specified
  by~\eqref{psi} when $\phi(x_1, x_2) \equiv p$.

  (a)~Consider a VERG $\bd G(3, \XX, \mu, \phi)$ representing an ERG
  ${\bd G(3, p)}$.  It can be shown easily that~$p$
  is an eigenvalue (say, $\lambda_1 = p$) with constant
  eigenfunction~${\bf 1}$.  
  [This  can be done by using the Cauchy--Schwarz inequality to prove that for any
  VERG with $n \geq 3$ we have the positive dependence
    \begin{equation}
      \label{posdep}
      \Pr\{1 \sim 2\mbox{\rm \ \ and\ \ }1 \sim 3\} \geq (\Pr\{1 \sim 2\})^2,
    \end{equation}
with equality if and only if the constant function~${\bf 1}$ is an eigenfunction
of~$T$ with eigenvalue $\Pr\{1 \sim 2\}$; moreover, we have equality
in~\eqref{posdep} for an ERG.
Cf.\ the proof of Theorem~\ref{thm:ERG-not-VRG}, where a similar argument is
used for $n\ge4$.]
  One then readily computes that the expected number of
  rooted cycles on three vertices is $6 \sum \lambda^3_i = 6 p^3$
  [this is Lemma~\ref{lem:rtdcyc} and~\eqref{ERGrc}, recalling that
  $n=3$] and similarly that the expected number of rooted edges is $6
  \lambda_1 = 6 p$ and the expected number of rooted paths on three
  vertices is $6 \lambda^2_1 = 6 p^2$.  So
\begin{equation}
\label{cubes}
\sum_{i \geq 2} \lambda^3_i = 0.
\end{equation} 
Conversely, suppose that a VERG $\bd G(3, \XX, \mu, \phi)$ has
eigenvalue $\lambda_1 = p$ with corresponding eigenfunction~$1$, and
that~\eqref{cubes} holds.  Then the expected counts of rooted edges,
rooted $3$-paths, and rooted $3$-cycles all agree with those for an
ERG ${\bd G(3, p)}$.  Since these three expected counts are easily
seen to characterize any isomorphism-invariant random graph model on
three vertices, the VERG represents the ERG ${\bd G(3, p)}$.

Summarizing, we see that a VERG $\bd G(3, \XX, \mu, \phi)$ represents
${\bd G(3, p)}$ if and only if $\lambda_1 = p$ with eigenfunction~$1$
and \eqref{cubes} holds.

In particular, one can take~$\mu$ to be the uniform distribution on
$\XX = [0, 1)$ and
$$
\phi(x_1, x_2) = g(x_1 \oplus x_2), \quad x_1, x_2 \in [0, 1),
$$
for any $g \geq 0$ satisfying $\int\!g(x)\,dx = p$.  It follows by
Lemma \ref{lem:unif} that then $\bd G(3, \XX, \mu, \phi)={\bd G(3,
  p)}$.  Alternatively, we can verify \eqref{cubes} by Fourier
analysis as follows.

Let $e_k(x)=e^{2\pi i k x}$. Then
\begin{equation*}
  (T e_k)(x) = \int_0^1 g(x\oplus y) e_k(y)\,dy 
  = \int_0^1 g(y) e_k(y-x)\,dy 
  = \hatg(-k) e_{-k}(x),
  \ \  k\in\mathbb Z.
\end{equation*}
For $k=0$, this says again that $e_0=1$ is an eigenfunction with
eigenvalue $\hatg(0)=p$.  For $k\ge1$, since
$\hatg(k)=\overline{\hatg(-k)}$, it follows that if $\omega_k
:=\hatg(k) / |\hatg(k)|$ (with $\omega_k := 1$ if this expression
would give $0 / 0$), then $\omega_ke_k\pm e_{-k}$ are eigenfunctions
with eigenvalues $\pm|\hatg(k)|$. Since $\{e_k\}$ is an orthonormal
basis, these eigenfunctions span a dense subspace of $L^2[0,1)$, so we
have found all eigenvalues, viz.\ $\lambda_1=p$ and $\pm|\hatg(k)|$,
$k=1,2,\dots$, and \eqref{cubes} follows.

(b)~The choice $g(x) = {\bf 1}(x \leq p)$ in~(a) was used
at~\eqref{psi} (when the VERG in question there is an ERG).  In this
case,
\begin{equation*}
  \hatg(k)
  =\int_0^p e^{-2\pi i k x}\,dx
  =\frac{1-e^{-2\pi i k p}}{2\pi i k}
\end{equation*}
and the multiset of eigenvalues can be listed as (changing the
numbering) $\{\lambda_j: j\in\mathbb Z\}$, where
\begin{equation*}
  \lambda_j := 
  \begin{cases}
    \frac{|1-e^{-2\pi i j p}|}{2\pi j}	
    =
    \frac{|\sin(\pi j p)|}{\pi j},
    &
    j \neq0,	
    \\
    p,& j=0.
  \end{cases}
\end{equation*}
\end{remark}

\section{Open problems}
\label{sect:open}

Call a VERG $\bd G(n, \XX, \mu, \phi)$ \emph{binary} if
$\Pr\{\phi(X_1, X_2) \in \{0, 1\}\} = 1$ where $X_1$ and $X_2$ are
independent draws from~$\mu$.  Since $\mu$-null sets do not matter,
this amounts to saying that~$\phi$ gives a representation of the
random graph as a VRG.  We will make use of the observation that
\begin{equation}
  \label{01}
  \mbox{~$\phi$ is binary\ \  if and only if\ \ }
  \EE [\phi(X_1, X_2) (1 - \phi(X_1, X_2))] = 0.
\end{equation} 

In Theorem~\ref{thm:VERG-is-VRG} we have seen that every VERG with $n
\leq 3$ is a VRG, but what is the situation when $n \geq 4$?

\begin{open}
  \label{open1}
  Is there any VRG with $n \geq 4$ that also has a non-binary VERG
  representation?
\end{open}

Theorem~\ref{thm:ERG-not-VRG} rules out constant-valued non-binary
VERG representations~$\phi$, and the main goal now is to see what
other VERGs we can rule out as VRGs.  In the following proposition,
$X_1$ and $X_2$ (respectively, $Y_1$ and $Y_2$) are independent draws
from~$\mu$ (respectively, $\nu$).

\begin{prop}
  If a VRG $\bd G(n, \YY, \nu, \psi)$ has a representation as a VERG
  $\bd G(n, \XX, \mu, \phi)$, then~$\phi$ is binary if and only if
  $\EE \psi^2(Y_1, Y_2) = \EE \phi^2(X_1, X_2)$.
\end{prop}

\begin{proof}
  Because $\bd G(n, \YY, \nu, \psi)$ and $\bd G(n, \XX, \mu, \phi)$
  represent the same random graph, we have
  $$
  \EE \psi(Y_1,
  Y_2) = \Pr\{1 \sim 2\} = \EE \phi(X_1, X_2).
  $$  
  Thus, by~\eqref{01}, $\phi$ is binary if and only if
  $$
  0 = \EE [\psi(Y_1, Y_2) (1 - \psi(Y_1, Y_2))] = \EE \psi(Y_1, Y_2) -
  \EE \psi^2(Y_1, Y_2)
  $$
  agrees with
  $$
  \EE [\phi(X_1, X_2) (1 - \phi(X_1, X_2))] = \EE \psi(Y_1, Y_2) - \EE  \phi^2(X_1, X_2),
  $$ 
  i.e.,\ if and only if $\EE \psi^2(Y_1, Y_2) = \EE \phi^2(X_1, X_2)$.
\end{proof}

The expression $\EE \phi^2(X_1, X_2)$ is the squared Hilbert--Schmidt
norm of the operator~$T$ defined at~\eqref{T} and equals the sum
$\sum_i \lambda^2_i$ of squared eigenvalues.  So the proposition has
the following corollary.

\begin{cor}
  \label{cor:squares}
  If a VRG $\bd G(n, \YY, \nu, \psi)$ has a representation as a VERG
  $\bd G(n, \XX, \mu, \phi)$, and if the respective multisets of
  nonzero squared eigenvalues of the integral operators associated
  with~$\psi$ and~$\phi$ are the same, then~$\phi$ is binary.~\qed
\end{cor}

\begin{open}
  \label{open2}
  Is there any VERG with $n \geq 4$ having two representations with
  distinct multisets of nonzero squared eigenvalues?
\end{open}

By Corollary~\ref{cor:squares}, a positive answer to Open
Problem~\ref{open1} would imply a positive answer to Open
Problem~\ref{open2}.

Our next result, Proposition~\ref{prop:rank1}, goes a step beyond
Theorem~\ref{thm:ERG-not-VRG}.  We say that~$\phi$ is of rank~$r$ when
the corresponding integral operator~\eqref{T} has exactly~$r$ nonzero
eigenvalues (counting multiplicities).  For~$\phi$ to be of rank at most~$1$ it is equivalent
that there exists $0 \leq g \leq 1$ ($\mu$-a.e.)\ such that (for
$\mu$-almost every~$x_1$ and~$x_2$)
\begin{equation}
  \label{rank1phi}
  \phi(x_1, x_2) = g(x_1) g(x_2).
\end{equation}

\begin{prop}
  \label{prop:rank1}
  For $n \geq 6$, no non-binary VERG $\bd G(n, \XX, \mu, \phi)$
  with~$\phi$ of rank at most~$1$ is a VRG.
\end{prop}

\begin{proof}
  Of course~$\phi$ cannot be both non-binary and of rank~$0$.  By
  Corollary~\ref{cor:squares} it suffices to show, as we will, that
\begin{quote}  
  $(*)$ \emph{any} VERG-representation $\bd G(n, \YY, \nu, \psi)$ of a VERG \\
  $\bd G(n, \XX, \mu, \phi)$ with $n \geq 6$ and~$\phi$ of rank~$1$
  must have the same single nonzero eigenvalue (without multiplicity).
\end{quote}
  Indeed, to prove~$(*)$, express~$\phi$ as at~\eqref{rank1phi} and let $\lambda_1,
  \lambda_2, \dots$ denote the eigenvalues corresponding to~$\psi$.
  By equating the two expressions for $\EE N_k$ obtained by applying
  Lemma~\ref{lem:rtdcyc} both to $\bd G(n, \XX, \mu, \phi)$ and to
  $\bd G(n, \YY, \nu, \psi)$, we find, with 
  $$
  c := \left[ \EE \phi^2(X_1, X_2) \right]^{1 / 2} > 0
  $$ 
  for shorthand, that
  \begin{equation}
  \label{rcidentity2}
  \sum_i \lambda^k_i = c^k, \qquad 3 \leq k \leq n.
  \end{equation}
  Applying~\eqref{rcidentity2} with $k = 4$ and $k = 6$, it follows
  from Lemma~\ref{lem:ineq} to follow (with $b_i := \lambda^4_i$ and
  $t = 3 / 2$) that~$\psi$ is of rank~$1$, with nonzero
  eigenvalue~$c$.
\end{proof}

The following lemma, used in the proof of
Proposition~\ref{prop:rank1}, is quite elementary and included for the reader's convenience.

\begin{lemma}
  \label{lem:ineq}
  If $b_1, b_2, \dots$ form a finite or infinite sequence of
  nonnegative numbers and $t \in (1, \infty)$, then
  $$
  \Big( \sum_i b_i \Big)^t \geq \sum_i b_i^t,
  $$
  with strict inequality if more than one $b_i$ is positive
and the right-hand sum is finite.
\end{lemma}

\begin{proof}
  The lemma follows readily in general from the special case of two
  $b$s, $b_1$ and $b_2$.  Since the case that $b_1 = 0$ is trivial, we
  may suppose that $b_1 > 0$.  Fix such a $b_1$, and consider the
  function
  $$
  f(b_2) := (b_1 + b_2)^t - b_1^t - b_2^t
  $$
  of $b_2 \geq 0$.  Then $f(0) = 0$ and
  $$
  f'(b_2) = t [(b_1 + b_2)^{t - 1} - b_2^{t - 1}] > 0.
  $$
  The result follows.
\end{proof}

With the hypothesis of Proposition~\ref{prop:rank1} strengthened to $n
\geq 8$, we can generalize that proposition substantially as follows.

\begin{prop}
  \label{prop:rankr}
  For $1 \leq r < \infty$ and $n \geq 4 (r + 1)$, no non-binary VERG
  $\bd G(n, \XX, \mu, \phi)$ with~$\phi$ of rank at most~$r$ is a VRG.
\end{prop}

It suffices to consider~$\phi$ of rank~$r$ exactly.  The strategy for
proving Proposition~\ref{prop:rankr} is essentially the same as for
Proposition~\ref{prop:rank1}: Under the stated conditions on~$n$
and~$r$, we will show that \emph{any} VERG-representation $\bd G(n,
\YY, \nu, \psi)$ of a VERG $\bd G(n, \XX, \mu, \phi)$ with~$\phi$ of
rank~$r$ must have the same finite multiset of nonzero squared
eigenvalues; application of Corollary~\ref{cor:squares} then completes
the proof.  The following two standard symmetric-function lemmas are
the basic tools we need; for completeness, we include their proofs.

\begin{lemma}
  \label{lem:symm1}
  Consider two summable sequences $a_1, a_2, \dots$ and $b_1, b_2,
  \dots$ of strictly positive numbers; each sequence may have either
  finite or infinite length.  For $1 \leq k < \infty$, define the
  elementary symmetric functions
  \begin{equation}
  \label{st}
  s_k := \sum_{i_1< i_2< \dots< i_k} a_{i_1} a_{i_2} \dots a_{i_k}, \qquad 
  t_k := \sum_{j_1< j_2< \dots< j_k} b_{j_1} b_{j_2} \dots b_{j_k}.
  \end{equation}
  For any $1 \leq K < \infty$, if $\sum_i a^k_i = \sum_j b^k_j$ for $k
  = 1, 2, \dots, K$, then {\rm (a)}~$s_k = t_k$ for $k = 1, 2, \dots,
  K$ and {\rm (b)}~the sequence ${\bd a}$ has length $\geq K$ if and
  only if the sequence ${\bd b}$ does.
\end{lemma}

\begin{proof}
  Clearly all the sums $\sum a^k_i$, $\sum b^k_j$, $s_k$, $t_k$ are
  finite, for any $k \geq 1$.  Using inclusion--exclusion, each $s_k$
  can be expressed as a finite linear combination of finite products
  of $\sum_i a^1_i$, $\sum_i a^2_i$, \dots $\sum_i a^k_i$.  (This is
  true when all indices~$i$ for $a_i$ are restricted to a finite
  range, and so also without such a restriction, by passage to a
  limit.)  Each $t_k$ can be expressed in just the same way, with the
  sums $\sum_j b^m_j$ substituting for the respective sums $\sum_i
  a^m_i$.  The assertion~(a) then follows; and since the
  sequence~${\bd a}$ has length $\geq K$ if and only if $s_K > 0$, and
  similarly for~${\bd b}$, assertion~(b) also follows.
\end{proof}

  \begin{lemma}
  \label{lem:symm2}
  Let $1 \leq K < \infty$, and let $a_1, \dots, a_K$ and $b_1, \dots,
  b_K$ be numbers.  If the sums~$s_k$ and~$t_k$ defined at~\eqref{st}
  satisfy $s_k = t_k$ for $k = 1, \dots, K$, then the multisets
  $\{a_1, \dots, a_K\}$ and $\{b_1, \dots, b_K\}$ are equal.
\end{lemma}

\begin{proof}
  We remark that the numbers $a_k$ and $b_k$ need not be positive, and
  may even be complex.  The result is obvious from the identity
  \begin{equation*}
    (z - a_1) \cdots (z - a_K) 
    = z^K - s_1 z^{K - 1} + s_2 z^{K - 2} + \cdots + (-1)^K s_K.~\qedhere
  \end{equation*}
\end{proof}
\smallskip

\begin{proof}[Proof of Proposition~\ref{prop:rankr}]
  Consider a VERG $\bd G(n, \XX, \mu, \phi)$ with~$\phi$ of rank~$r$,
  and let $M = \{\lambda^2_1, \lambda^2_2, \dots, \lambda^2_r\}$ be
  its multiset of nonzero squared eigenvalues.  Suppose that the same
  random graph can also be represented as the VERG $\bd G(n, \YY, \nu,
  \psi)$, and let the finite or infinite multiset $\widetilde{M} :=
  \{\tilde{\lambda}^2_1, \tilde{\lambda}^2_2, \dots \}$ be the
  multiset of nonzero squared eigenvalues for~$\psi$.  As discussed
  immediately following the statement of the proposition, it suffices
  to show that the multisets~$M$ and~$\widetilde{M}$ are equal.

   Let $a_i := \lambda^4_i$ and $b_j := \tilde{\lambda}^4_j$.  Applying
  Lemma~\ref{lem:rtdcyc} with $k = 4, 8, \dots, 4 (r + 1)$, we see
  that the hypotheses of Lemma~\ref{lem:symm1} are satisfied for $K =
  r$ and for $K = r + 1$.  Therefore, $\widetilde{M}$ has size~$r$ and
  the sums~\eqref{st} satisfy $s_k = t_k$ for $k = 1, 2, \dots, r$.
  By Lemma~\ref{lem:symm2}, the two multisets are equal.
\end{proof}
\medskip

{\bf Acknowledgment.}\ \ The authors thank an anonymous reviewer
who provided us with helpful feedback on an earlier version of this paper.

\bibliographystyle{plain} 
\bibliography{paper20101012}
\end{document}